\documentclass[12pt]{amsart}
\usepackage{fancyhdr}
\usepackage{amssymb,latexsym,amscd}
\usepackage[cp850]{inputenc}
\newcommand{\cO}{\mathcal{O}}

\newcommand{\lra}{\longrightarrow}
\newcommand{\hra}{\hookrightarrow}
\newcommand{\ra}{\rightarrow}

\newcommand{\ms}{\mapsto}
\newcommand{\tm}{\subset}

\newcommand{\CC}{\mathbb C}
\newcommand{\RR}{\mathbb R}

\newcommand{\PP}{\mathbb P}

\theoremstyle{plain}

\begin{document}
\title[Coherent Systems of genus 0]{Coherent Systems of genus 0}
\author{H. Lange}
\author{P.E. Newstead}
\address{Mathematisches Institut\\
              Universit\"at Erlangen-N\"urnberg\\
              Bismarckstra\ss e $1\frac{ 1}{2}$\\
              D-$91054$ Erlangen\\
              Germany}
\email{lange@mi.uni-erlangen.de}
\address{Department of Mathematical Sciences\\University of Liverpool\\
           Peach Street\\Liverpool L69 7ZL\\
           U.K.}
\email{newstead@liv.ac.uk}
\thanks{Supported by DFG Contracts Ba 423/8-1 and EPSRC Grant no. GR/S48974. Both authors are members
of the research group VBAC (Vector Bundles on Algebraic Curves), which is
partially supported by EAGER (EC FP5 Contract no. HPRN-CT-2000-00099) and
by EDGE (EC FP5 Contract no. HPRN-CT-2000-00101). The second author would also like to thank the Royal Society
and CSIC, Madrid, for supporting a visit to CSIC during which this work was completed.}
\begin{abstract}
In this paper we begin the classification of coherent systems $(E,V)$ on the projective line which are stable with respect to
some value of a parameter $\alpha$. In particular we show that the moduli spaces, if non-empty, are always smooth and
irreducible of the expected dimension. We obtain necessary conditions for non-emptiness and, when $\dim V=1$ or $2$, 
we determine these conditions precisely. We also obtain partial results in some other cases. 
\end{abstract}
\maketitle


\section{Introduction}

A {\it coherent system of type $(n,d,k)$} on a smooth projective curve $C$ over an algebraically closed field is by 
definition a pair $(E,V)$ with 
a vector bundle $E$ of rank $n$ and degree $d$ over $C$ and a vector subspace $V \tm H^0(E)$ of dimension $k$. 
For any real number $\alpha$, the {\it $\alpha$-slope} of a coherent system $(E,V)$ of type $(n,d,k)$ is defined by
$$
\mu_{\alpha}(E,V) := \frac{d}{n} + \alpha \frac{k}{n}
$$
A {\it coherent subsystem} of $(E,V)$ is a coherent system $(E',V')$ such that $E'$ is a subbundle of $E$ and $V' \subset V \cap H^0(E')$.
A coherent system $(E,V)$ is called 
{\it $\alpha$-stable} ($\alpha$-semistable) if
$$
\mu_{\alpha}(E',V') < \mu_{\alpha}(E,V) \ \ (\mu_{\alpha}(E',V') \le \mu_{\alpha}(E,V))
$$
for every proper coherent subsystem $(E',V')$ of $(E,V)$.
The $\alpha$-stable coherent systems of type $(n,d,k)$ on $C$ form a quasiprojective moduli space which we denote by $G(\alpha;n,d,k)$.
These spaces have attracted a good deal of attention in the last few years because of applications to the study of
the moduli of vector bundles and also because they arise in gauge theory in connection with a generalisation of the vortex
equations \cite{5}.  At least in the case $k=1$, coherent systems can also be interpreted in terms of bundles 
on $C\times\PP^1$ which are equivariant for the action of $SL(2)$ \cite{6}.

A systematic study has been started in \cite{2}, to which we refer for general information on coherent systems. The general
theory is true for curves of all genera, but most of the detailed results require $g\ge2$. In the current 
paper, on the other hand, we study the moduli spaces $G(\alpha;n,d,k)$ on the projective line $\PP^1$. 
In this case there are no stable vector bundles
of rank $\ge2$, but the spaces $G(\alpha;n,d,k)$ may be nonempty, as we shall see, and often have different features from
those of higher genus. 

After some preliminaries, we prove in section 3 that all non-empty $G(\alpha;n,d,k)$ are smooth and irreducible of the
expected dimension; moreover, for a general $(E,V)\in G(\alpha;n,d,k)$, $E$ is of generic splitting type (see section 2 for
definitions). We give standard forms for the general $(E,V)$ in all cases. In section 4, we obtain conditions on $\alpha$ 
which are necessary for the existence of $\alpha$-stable bundles. Then, in section 5, we obtain precise conditions for
the non-emptiness of $G(\alpha;n,d,k)$ when $k=1$ and when $k=2$. Finally, we obtain partial results for some higher
values of $k$ in section 6.

We suppose throughout that $n\ge2$.

\section{Preliminaries}

We recall first that necessary conditions for $G(\alpha;n,d,k)$ to be nonempty are $\alpha>0$ and $d > 0$ 
\cite[Lemmas 4.2 and 4.3]{2}.  In genus $0$, the {\it Brill-Noether number} $\beta(n,d,k)$ \cite[Definition 2.7]{2} of 
coherent systems of type $(n,d,k)$ is given by
$$
\beta(n,d,k) = -n^2 + 1 - k(k-d-n).
$$
Note that, for any component $U$ of $G(\alpha;n,d,k)$, we have (see \cite[Corollary 3.6]{2})
$$
\dim U \geq \beta(n,d,k).  
$$
Moreover $G(\alpha;n,d,k)$ is smooth of dimension $\beta(n,d,k)$ at $(E,V)$ if and only if the {\it Petri map}
$$
V\otimes H^0(E^*\otimes K)\ra H^0(E\otimes E^*\otimes K)
$$
(given by multiplication of sections, and where $K$ denotes the canonical bundle) is injective (see \cite[Proposition 3.10]{2}).
Of course in our case $K\simeq\cO(-2)$.\\

We need also to recall some facts about vector bundles on $\PP^1$. Every vector bundle $E$ on $\mathbb P^1$ can be written 
uniquely as 
$$
E \simeq \bigoplus_{i=1}^{n}\cO(a_1)\ \ \mbox{with}\ \ a_1 \geq a_2 \geq \cdots \geq a_n.
$$
Given $n$ and $d$, we say that a bundle $E$ of rank $n$ and degree $d$ is of {\it generic splitting type} if
$a_1\le a_n+1$, in other words if $E$ can be written as
$$
E \simeq \cO(a+1)^s \oplus \cO(a)^{n-s},
$$
where $a$ and $s$ are defined by $d=an+s$ with $0 \leq s < n$. In particular the isomorphism class of such a bundle $E$ 
is determined by $n$ and $d$. Moreover $E$ is of generic splitting type if and only if $h^1(End(E)) = 0$. From this it follows 
that, in a family of bundles, those of generic splitting type form an open subset.\\

\section{Irreducibility of moduli spaces}

Fix a pair $(n,d)$ with $n\ge2$. In this section we show that the moduli spaces $G(\alpha;n,d,k)$ are irreducible
for any  $k$ whenever they are non-empty and study the generic elements $(E,V)$ of $G(\alpha;n,d,k)$.\\

{\bf Lemma 3.1:} {\it Suppose $k > 0$ and $(E,V)$ is $\alpha$-stable for some $\alpha$. Then
$$
E \simeq \bigoplus_{i=1}^{n} \cO(a_i)
$$
with all $a_i \geq 1$.}\\

{\it Proof:} We require to prove that $E$ has no direct factor $\cO(b)$ with $b \leq 0$.  If $b < 0$, 
we would have $(\cO(b),0)$ as a direct factor
of the coherent system $(E,V)$ contradicting $\alpha$-stability for all $\alpha$. If $b=0$ and the induced map 
$V \ra H^0(\cO(b))$ is $0$, the same applies. Otherwise $(\cO(b),H^0(\cO(b)))\simeq(\cO,H^0(\cO))$ is a direct factor of $(E,V)$, again
contradicting $\alpha$-stability for all $\alpha$.          \hfill $\square$\\

{\bf Theorem 3.2:} {\it Suppose $k>0$ and $G(\alpha;n,d,k)$ is non-empty. Then $G(\alpha;n,d,k)$ is 
smooth and irreducible and has dimension $\beta(n,d,k)$. Moreover,  
for a general $(E,V) \in G(\alpha;n,d,k)$, $E$ is of generic splitting type.}\\

{\it Proof}: If $(E,V)$ is $\alpha$-stable, it follows from Lemma 3.1 that $H^0(E^* \otimes K) = 0$, where $K$
denotes the canonical line bundle on $\mathbb P^1$. Hence the Petri map $V \otimes H^0(E^* \otimes K) \to
H^0(E \otimes E^* \otimes K)$ is injective. So $G(\alpha;n,d,k)$ is smooth and has dimension 
$\beta(n,d,k)$.
By Lemma 3.1 there can be only finitely many $E$ such that $(E,V) \in G(\alpha;n,d,k)$ for some 
$V \subset H^0(E)$. For any such $E$, consider the Grassmannian $Gr(k,H^0(E))$. Define
an open subset  $U$ of  $Gr(k,H^0(E))$ by 
$$
U = \{ V \,\,|\,\, (E,V)\,\,\, \mbox{is}\,\,\, \alpha\mbox{-stable}\}.
$$ The canonical morphism $U \to G(\alpha;n,d,k)$
has fibres which are orbits for the natural action of $\mbox{Aut}(E)$ on $U$. Since $(E,V)$ is $\alpha$-stable, the stabilizers
for this action consist of scalar multiples of the identity (see \cite[Proposition 2.2]{2}). 
So the image $U_1$ of $U$ in $G(\alpha;n,d,k)$ has 
dimension 
$$
\begin{array}{ll}
\dim U_1 &= \dim Gr(k,H^0(E)) - \dim \mbox{Aut}(E) + 1\\
&= k(d + n -k) - n^2 - h^1(End(E)) + 1\\
&= \beta(n,d,k) - h^1(End(E)).
\end{array}
$$
Hence $\overline{U_1}$ cannot contain a component of $G(\alpha;n,d,k)$ unless $h^1(End(E)) = 0$, i.e. $E$ is 
of generic splitting type. This determines $E$, thus giving the irreducibility of $G(\alpha;n,d,k)$. 
The other assertions follow immediately. \hfill $\square$\\

{\bf Corollary 3.3:} {\it If $G(\alpha;n,d,k) \not= \emptyset $, then $ d \geq \frac{1}{k}(n^2 - 1) - (n-k)$}.\\

{\it Proof}: The stated inequality for $d$ is equivalent to $\beta(n,d,k)\ge0$.\hfill$\square$\\

{\bf Corollary 3.4:} {\it For any $(n,d,k)$, the set $\{\alpha\,\,|\,\,G(\alpha;n,d,k)\ne\emptyset\}$ is an open interval in
$\RR$ (possibly empty or semi-infinite).}\\

{\it Proof}: By the theorem and the fact that, as $\alpha$ varies, the $\alpha$-stability condition can only change as $\alpha$
passes through one of a finite number of critical values \cite[Propositions 4.2 and 4.6]{2}, the general element of
$G(\alpha;n,d,k)$ is independent of $\alpha$. The result now follows from \cite[Lemma 3.14]{7}.\hfill$\square$\\

{\bf Proposition 3.5:} 
{\it Let $0 < k < n$ and suppose $G(\alpha;n,d,k)$ is non-empty. Then for a general $(E,V) \in G(\alpha;n,d,k)$ we have an exact sequence
$$
0 \lra \cO^k \lra E \lra G \lra 0  \eqno(1)
$$
where $V = H^0(\cO^k) \subset H^0(E)$ and $G$ is a vector bundle. Moreover
$$
E \simeq \cO(a)^{n-t} \oplus \cO(a-1)^t \quad \mbox{and} \quad   
G \simeq \cO(a+l+1)^m \oplus \cO(a+l)^{n-k-m}   \eqno(2)
$$
where the integers $a,t,l$ and $m$ are defined by
$$
d=na-t \,\, \mbox{with} \,\, 0 \leq t < n \,\, \mbox{and} \,\, ka-t = l(n-k) + m \,\, \mbox{with} \,\,
0 \leq m < n-k.   \eqno(3)
$$}

{\it Proof}: By Theorem 3.2, we may assume that $E$ is of generic splitting type. Moreover, by Lemma 3.1,
$E$ is generated by its global sections. Hence by a standard result (see \cite[Theorem 2]{1}), we obtain the
exact sequence (1) for a generic subspace $V$ of $H^0(E)$ of dimension $k$. 
It remains to show that $G$ has the stated form.
To see this, we consider for any $G$ the subscheme $U(G)$ of $G(\alpha;n,d,k)$ given by the exact
sequences 
$$
0 \lra \cO^k \lra E' \lra G \lra 0.
$$
Note that $h^0(G^*) = 0$, for otherwise $(E',V)$ would have an endomorphism not equal to a scalar multiple of 
the identity, contradicting the $\alpha$-stability of $(E',V)$. Hence
$$
\begin{array}{ll}
\dim U(G) &= \dim \mbox{Ext}^1(G,\cO^k) - \dim \mbox{Aut}(\cO^k) -  \dim \mbox{Aut}(G) + 1\\
&= k \cdot h^1(G^*) - k^2 - (n-k)^2 - h^1(End(G)) + 1\\
&= k(d-n+k) - 2k^2 - n^2 + 2nk - h^1(End(G)) + 1\\
&= \beta(n,d,k) - h^1(End(G)).
\end{array}
$$
So, for $(E,V)$ general, we must have $h^1(End(G)) = 0$, i.e. $G$ is of generic splitting type. This completes the 
proof.  \hfill $\square$\\

{\bf Proposition 3.6:} {\it Suppose $G(\alpha;n,d,n)$ is non-empty. Then for a general $(E,V) \in G(\alpha;n,d,n)$
we have an exact sequence
$$
0 \lra \cO^n \lra E \lra T \lra 0 \eqno(4)
$$
where $V = H^0(\cO^n) \subset H^0(E)$ and $T$ is a torsion sheaf of length $d$.}\\

{\it Proof:} According to Lemma 3.1 $E$ is generated by its global sections.
If we fix a point $p$ of ${\mathbb P}^1$ and choose $V$ generically, the evaluation map of $V$ at $p$ is an 
isomorphism. Now the choice of a basis of $V$ defines an exact sequence (4).  \hfill $\square$\\

{\bf Proposition 3.7:} {\it Suppose $k > n$ and $G(\alpha;n,d,k)$ is non-empty. For the general $(E,V) \in 
G(\alpha;n,d,k)$ we have an exact sequence 
$$
0 \lra H \lra \cO^k \lra E \lra 0   \eqno(5)
$$
where $H^0(\cO^k) \to H^0(E)$ is an isomorphism onto $V$}.\\

{\it Proof:} According to Lemma 3.1 $E$ is generated by its global sections. Then by a standard result 
(easily obtained from \cite[Theorem 2]{1}), the general subspace of $H^0(E)$ of dimension $k$ defines 
a surjection $\cO^k \to E$. \hfill $\square$\\

\section{Bounds for $\alpha$}

We recall that in general non-empty moduli spaces can exist only if $\alpha>0$. Moreover, for $0<k<n$, 
we must have also $\alpha<\frac{d}{n-k}$ \cite[Lemma 4.2]{2}.
In this section we show that, in our case, the assumption $G(\alpha;n,d,k) \not= \emptyset$ can imply stronger upper and lower 
bounds on  $\alpha$.\\

{\bf Proposition 4.1:} {\it Suppose $G(\alpha;n,d,k)$ is non-empty and let the numbers $a$ and $t$ be defined  
by $d = na - t$ with $0 \leq t < n$. Then
$$
\alpha > \frac{t}{k}.  \eqno(6)
$$
Moreover, if $E$ has generic splitting type and (6) holds, then every coherent subsystem of $(E,V)$ 
of type $(r,d',0)$ satisfies the $\alpha$-stability condition.}\\ 

{\it Proof}: Let $(E,V) \in G(\alpha;n,d,k)$. We may assume that $E$ has generic splitting type 
$E \simeq \cO(a)^{n-t} \oplus \cO(a-1)^t$. Then $(\cO(a),0)$ is a proper coherent subsystem of $(E,V)$ and the $\alpha$-stability condition for it gives (6). For the last fact, we note that every subbundle $F$ of $E$ 
has slope
$\leq a$ and apply the $\alpha$-stability condition to $(F,0)$.  \hfill $\square$\\

{\bf Proposition 4.2:} {\it Suppose $0 < k < n$ and $G(\alpha;n,d,k)$ is non-empty and let the integers $a,t,l$ 
and $m$ be as in (3). Then
$$
\alpha < \frac{d}{n-k} - \frac{mn}{k(n-k)}.  \eqno(7)
$$
Moreover, if $E$ and $G$ are as in (2) and (7) holds, then every coherent subsystem of $(E,V)$ of type 
$(r,d',k)$ satisfies the $\alpha$-stability condition.}\\

{\it Proof}: Let $(E,V) \in G(\alpha;n,d,k)$. By Proposition 3.5, we may assume that there is an 
exact sequence (1) with $E$ and $G$ as in (2).
Let $(F,W)$ be a coherent subsystem of rank $r$ with $\dim W = k$.
Then we have a diagram
$$
\begin{array}{ccccccccc}
&&0&&0&&&&\\
& & \downarrow & & \downarrow & & & & \\
0& \ra & \cO^k & \stackrel{\sim}{\ra} & \cO^k & \ra & 0 & &\\
& & \downarrow & & \downarrow & & \downarrow & & \\ 
0 & \ra & F & \ra & E & \ra & H & \ra & 0 \\
& & \downarrow & & \downarrow & & || & & \\
0 & \ra & K & \ra & G & \ra & H & \ra & 0 \\
& & \downarrow & & \downarrow & & \downarrow & & \\
&&0&&0&&0&&\\
\end{array}
$$
in which $K$ and $H$ are vector bundles.

If $r \geq k + m$, then $\deg F = \deg K \leq (r-k)(a+l+1) -(r-k-m)$
where equality can be attained. The stability of $(E,V)$ implies that
$$
\mu_{\alpha}(F,W) = \frac{1}{r}((r-k)(a+l+1) - (r-k-m)) + \frac{k\alpha}{r} < \mu_{\alpha}(E,V) = \frac{d}{n} + \frac{k\alpha}{n},
$$
which is equivalent to
$$
\begin{array}{ll}
 k (n-r)  \alpha &< rd - (a+l+1)(r-k)n + (r-k-m)n\\
& = rd-(a+l)(r-k)n-mn.
\end{array}
$$
Multiplying this inequality by $n-k$ and using the equation 
$$
(a+l)(n-k) = a(n-k)+ka - t - m=d-m,
$$ 
this is equivalent to
$$
\begin{array}{ll}
k(n-r)(n-k) \alpha &< r(n-k)d -(r-k)n(d-m) -mn(n-k)\\
& = k(n-r)d - mn(n-r).
\end{array}
$$
Dividing by $k(n-r)(n-k)$, this gives (7).
This proves the first part of the proposition, and also the second part under the 
assumption $r \geq k+m$.\\

If $r < k+m$, then $\deg F = \deg K \leq (r-k)(a+l+1)$ where again equality can be attained. The stability of 
$(E,V)$
implies that
$$
\mu_{\alpha}(F,W) = \frac{r-k}{r}(a+l+1) + \frac{k\alpha}{r} < \mu_{\alpha}(E,V) = \frac{d}{n} + \frac{k\alpha}{n},
$$
which is equivalent to
$$
k(n-r)\alpha < dr - (a+l+1)(r-k)n
$$
The same substitution as above now gives
$$
k(n-r)(n-k) \alpha < k(n-r)d - n(r-k)(n-m-k)
$$
and hence to
$$
\alpha < \frac{d}{n-k} - \frac{n(r-k)(n-m-k)}{k(n-k)(n-r)}.  \eqno(8)
$$
It remains to show that for $k \leq r < k+m$ the right hand side of (7) is less or equal to the right hand side
of (8). But this is an immediate computation.   \hfill $\square$\\

{\bf Remark 4.3:}  Suppose $0 < k < n$. Then the existence of $\alpha$ for which (6) and (7) both hold is 
equivalent to $(n-k)t < kd - mn$. Substituting for $t$ and $m$ from (3), this is equivalent to $l>0$.\\

{\bf Proposition 4.4:} {\it Let $l$ be as in (3) and suppose $l > 0$. Then there exists a sequence (1) with
$E$ and $G$ as in (2).}\\

{\it Proof}: Suppose $l>0$. Note that (3) implies $a>0$. Now a theorem of Shatz  \cite[Theorem 1]{3} ensures the existence of an
exact sequence $0 \ra \cO^k \ra E \ra G \ra 0$ if the following two conditions are satisfied  for the 
Harder-Narasimhan polygons 
$HNP$ of $E$ and $F := G \oplus \cO^k$
\begin{itemize}
\item $HNP(F) \geq HNP(E)$, 
\item if $E = \oplus_{i=1}^n \cO(a_i)$ with $a_1\ge\ldots\ge a_n$ and $F = \oplus_{i=1}^n \cO(b_i)$ with 
$b_1\ge\ldots\ge b_n$, then 
$ b_i > a_i$ exactly for $i \leq n-k$.
\end{itemize}
These conditions are certainly fulfilled for $E$ and $G$ as in (2)
if $l>0$. \hfill $\square$\\

\section{Coherent systems with $k= 1$ and $k=2$}

In this section we determine the $\alpha$ for which $G(\alpha;n,d,k)$ is non-empty in the cases $k=1$ and $k=2$.
For $k=1$ we obtain as an immediate consequence of Propositions 3.1, 3.2 and 3.4:\\

{\bf Theorem 5.1:} {\it Suppose $n \geq 2$.  Then {\it $G(\alpha;n,d,1) \not= \emptyset$ if and only if} 
$$t < \alpha < \frac{d}{n-1} - \frac{mn}{n-1}.$$}\\

{\bf Remark 5.2:} For $k=1$, (3) can be solved to give
$$
d = n(n-1)l + mn + t(n-1).       \eqno(9)
$$
So $l>0$ implies $d \geq n^2 - n$ (compare Corollary 3.3), but in general $d$ does not take all values
$\geq n^2 - n$. The precise values for which $G(\alpha,n,d,1) \not= \emptyset$ for some $\alpha$ are given
by (9) with $l>0, \,\,0 \leq t<n$ and $0 \leq m < n-1$.\\ 

Before we can prove the analogous statement in the case $k=2$, we need an auxiliary result. Let 
$(E,V)$ be a coherent system as in Proposition 3.5 with $k=2$ and let 
$\beta: \cO^2 \ra E \simeq \cO(a)^{n-t} \oplus \cO(a-1)^t$ denote 
the associated embedding. For $t \geq 1$ consider the invariant $\delta$ of the coherent system $(E,V)$ defined as follows
$$
\delta(E,V) = \left\{ \begin{array}{l} \mbox{minimal rank of a direct factor of} \,\, \cO(a-1)^t \\ 
\mbox{containing the image of some} \,\, \cO \subset \cO^2 \,\, \mbox{under}\\ \mbox{the composed map} \,\,
\cO \hra \cO^2 \stackrel{\beta}{\ra} E \ra \cO(a-1)^t. \end{array} \right.
$$\\

{\bf Lemma 5.3:} {\it Suppose $a \geq 1$ and $t \geq 1$. Then the general coherent system $(E,V)$ of type
$(n,d,2)$ as in Proposition 3.5 satisfies}
$$\delta(E,V) = \left\{ \begin{array}{lll} t & & a \geq t+1,\\
                                                      t-1& if & a=t,\\
                                                       a && 1 \leq a < t.
                             \end{array} \right.
$$\\

{\it Proof}: The map $\beta$ is given by a matrix of the form 
$$
M = \left( \begin{matrix}
f_1&\cdots&f_{n-t}&g_1&\cdots&g_t \cr
            f'_1&\cdots&f'_{n-t}&g'_1&\cdots&g'_t \cr
\end{matrix}  \right)^{tr}          
$$
with binary forms $f_i, f'_i$ of degree $a$ and $g_i, g'_i$ of degree $a-1$ (where $tr$ denotes the transpose 
of the matrix). The composition $\cO \lra \cO(a-1)^t$ is then given
by the column vector $v = (bg_1 + cg'_1, \ldots, bg_t + cg'_t)^{tr}$ for some constants $b,c \in \CC$ with $(b,c) \not= (0,0)$. 
Now by definition of $\delta(E,V)$ we have
$$
t-\delta(E,V) = \max_{A \in GL(t,\mathbb C), \,\,(0,0) \not= (b,c) \in \mathbb C^2} \,\, \{\mbox{number of zero entries in} \,\, Av \}.
$$
But this equals the maximum number of linearly independent vectors 
$(\lambda_1, \ldots , \lambda_t) \in \mathbb C^t$ such that
$$
(b\lambda_1 g_1 + c\lambda_1 g'_1) + \cdots + (b \lambda_t g_t + c\lambda_t g'_t) = 0   \eqno(10)
$$
the maximum to be taken over all   $(b,c) \in \mathbb C^2$, $(b,c)\ne(0,0)$. Now consider the Segre embedding
$$
i: \mathbb P^1 \times \mathbb P^{t-1} \lra \mathbb P^{2t-1}
$$
$$
i((b:c),(\lambda_1: \ldots: \lambda_t)) = (b\lambda_1:c\lambda_1: \ldots : b\lambda_t:c\lambda_t).
$$
If we let $W$ denote the kernel of the linear map $\mathbb C^{2t} \lra H^0(\cO(a-1))$ given by 
$$
(\mu_1,\nu_1, \ldots, \mu_t,\nu_t) \ms \mu_1g_1 + \nu_1g'_1 + \cdots + \mu_tg_t + \nu_tg'_t ,
$$
then we have
$$
t - \delta(E,V) - 1 \leq \dim( \mathbb P(W) \cap i(\mathbb P^1 \times \mathbb P^{t-1}))  \eqno(11)
$$
Note that, for a general choice of $g_1,g'_1, \ldots ,  g_t,g'_t$ we have $ \dim W = \max\{2t-a,0\}$,
and that by varying $(E,V)$ generically we can choose $W$ arbitrarily.\\

If $a \geq t+1$, we have 
$$
2t - a - 1 + t < 2t - 1,
$$ 
and we can choose $W$ so that $\mathbb P(W) \cap i(\mathbb P^1 \times \mathbb P^{t-1}) =\emptyset$.
Hence (11) implies that $\delta(E,V) \geq t$, and so $\delta(E,V) = t$.\\ 

If $a=t$, then we can choose $W$ so that $\mathbb P(W) \cap i(\mathbb P^1 \times \mathbb P^{t-1})$
is finite and non-empty and thus by (11)  $\delta(E,V) \geq t-1$. On the other hand, any point of this intersection 
gives a solution of (10). So $t- \delta(E,V) \geq 1$, which gives the assertion in this case.\\

Finally suppose $1 \leq a < t$. Then we can choose $W$ so that 
$$
\dim(\mathbb P(W) \cap i(\mathbb P^1 \times \mathbb P^{t-1})) = (2t - a -1) +t-(2t-1) = t-a.
$$
Moreover this intersection is irreducible by Bertini. Hence the maximal dimension of a linear space contained
in $\mathbb P(W) \cap i(\mathbb P^1 \times \mathbb P^{t-1})$ is exactly $t-a-1$. This proves that $t-\delta(E,V) -1 = t-a-1$
which completes the proof of the lemma.  \hfill $\square$\\

{\bf Theorem 5.4:} {\it Suppose $n\ge3$.  Then
$G(\alpha;n,d,2) \not= \emptyset$ if and only $\alpha$ satisfies the condition
$$
\frac{t}{2} < \alpha < \frac{d}{n-2} - \frac{mn}{2(n-2)},   \eqno(12)
$$ 
the Brill-Noether condition
$$
d\ge\frac12n(n-2)+\frac32\eqno(13)
$$
and $(n,d)\ne(4,6)$.}\\

{\it Proof}: If $G(\alpha;n,d,2) \not= \emptyset$ then (12) holds by Propositions 4.1 and 4.2
and (13) by Corollary 3.3. If $(n,d)=(4,6)$ and $(E,V)$ is general, we have $E\simeq\cO(2)^2\oplus\cO(1)^2$ and it is easy to see 
that $(E,V)$ possesses a subsystem $(F,W)$ with $F\simeq\cO(2)\oplus\cO(1)$ and $\dim W=1$. This contradicts the
$\alpha$-stability condition for all $\alpha$.

For the converse, note that we can assume $l\ge1$ by Remark 4.3. Moreover, by Proposition 4.4, there exists
an exact sequence (1) with $E$ and $G$ as in (2) and $a$, $t$, $l$ and $m$ defined by (3). It suffices to prove
that the general $(E,V)$ of this form is $\alpha$-stable if (12) and (13) hold and $(n,d)\ne(4,6)$.

According to the last sentences of Propositions 4.1 and 4.2 $(E,V)$ satisfies the $\alpha$-stability condition for every coherent
subsystem of the form $(F,0)$ and for those of the form $(F,W)$ with $\dim W = 2$.
It remains to consider the coherent subsystems $(F,W)$ with $\mbox{rk} F = r$ and $\dim W = 1$. We can also suppose
$r\ge2$, since otherwise $(F,W)=(\cO,H^0(\cO))$, which certainly satisfies the $\alpha$-stability condition when (12) holds. 

In completing the proof, we distinguish
several cases\\

{\bf Case 1:} $r > \frac{n}{2}.$\\

Suppose first $t=0$. Then $E \simeq \cO(a)^n$ and any subbundle $F$ of $E$ is of slope $\leq a$. So $F$ satisfies the $\alpha$-stability 
condition if $\frac{\alpha}{r} < \frac{2\alpha}{n}$ which is equivalent to $2r >n$.

Hence we may assume $t \geq 1$. 
In this case we consider $H$, the subbundle image of $F$ in $\cO(a-1)^t$, and write $s= \mbox{rk} H$. According
to the definition of $\delta = \delta(E,V)$,
$$
h^0(H^*(a-1)) \geq \delta.
$$
On the other hand by Riemann-Roch $h^0(H^*(a-1)) = -\deg H +sa$, implying
that $\deg H \leq sa - \delta$ and thus
$$
\deg F \leq (r-s)a +sa -\delta = ra - \delta.
$$
Hence if we have 
$a - \frac{\delta}{r} + \frac{\alpha}{r} < \frac{d}{n} + \frac{2\alpha}{n} = a - \frac{t}{n} + \frac{2\alpha}{n}$
or equivalently
$$
(2r - n)\alpha > rt -n\delta, \eqno(14)
$$
the $\alpha$-stability condition for $(F,W)$ is certainly satisfied.

But by (12) $\,\,\, \alpha > \frac{t}{2} \,\,\,$ and thus
$$
(2r-n)\alpha > rt - n\frac{t}{2}.
$$
Thus it is enough to prove that $2\delta\ge t$. For this we apply Lemma 5.3. Note first that, since $l \geq 1$ 
and $n\ge3$,
$$
2a-t = l(n-2) + m \geq n-2 > 0.
$$
Hence
$$
a > \frac{t}{2}.  \eqno(15)
$$
Now if $a \geq t+1$, then $\delta = t$ and so $2\delta>t$.
If $a \leq t-1$, then $\delta = a$ and thus $2\delta>t$ by (15). Finally if $a = t$, then
$\delta = t-1$ and thus 
$$
2\delta=2t-2\ge t
$$ 
provided that $t=a\ge2$. However, if $t=a=1$, then (3) gives $l(n-2)+m=1$. For $l\ge1$, this can happen only if
$n=3$; but then $d=2$, which contradicts (13).\\

{\bf Case 2:} $r = \frac{n}{2}.$\\

Suppose first $t \geq 1$. As in Case 1 we have only to verify (14) which in 
this case says
$$
 2\delta > t.
$$
This follows by the same argument as in Case 1 except when $a=t=2$. But then (3) implies that $2 = l(n-2)+m$, which, 
since $l \geq 1$, can be solved only for $n=4$. But then 
$d = na -t = 6$, so we are in the exceptional case of the theorem.

Now let $t=0$, i.e $E \simeq \cO(a)^{n}$. In this case we use a modified version of the invariant $\delta$ of $(E,V)$. 
In fact, define $\delta' = \delta'(E,V)$  by
$$
 \delta'(E,V) = \left\{ \begin{array}{l} \mbox{minimal rank of a direct factor of} \,\, E \simeq \cO(a)^n\\ 
\mbox{containing the image of some} \,\, \cO \subset \cO^2 \,\, \mbox{under}\\ \mbox{the composed map} \,\,
\cO \hra \cO^2 {\hra} E. \end{array} \right.
$$\\
Lemma 5.3 also applies in this case and gives 
$$
\delta'(E,V) = \left\{ \begin{array}{lll} 2r & & a \geq 2r,\\
                                                      2r-1& if & a=2r-1,\\
                                                       a+1 && 1 \leq a \leq 2r-2.
                             \end{array} \right.   \eqno(16)
$$
According to the definition of $\delta'$
$$
h^0(F^*(a)) \geq \delta'.
$$
On the other hand, writing $e = \deg F$, we have by Riemann-Roch $h^0(F^*(a)) = -e + (a+1)r$, implying that
$$
e \leq (a+1)r - \delta'.   \eqno(17)
$$
We have to verify the $\alpha$-stability condition which is equivalent to $e < ar$.\\

If $a \geq 2r$, (16) and (17) give $e \leq (a-1)r < ar$. If $ a =2r-1$, we get
$e \leq (a-1)r+1 < ar$ since $r\ge2$. Finally, if $1 \leq a \leq 2r-2$, we get 
$$e \leq (a+1)r - (a+1) = ar -a +r -1.
$$
So it remains to show that $a \geq r$. But, by (13),
$$ \,\, 2ra = d \geq \frac{1}{2}2r(2r-2) + \frac{3}{2},
$$ 
which implies the assertion.\\

{\bf Case 3:} $r < \frac{n}{2}$.\\

In this case we consider the image $K$ of $F$ in $G$, which we may assume to be a subsheaf of $G$ 
of degree $e$ and rank $r-1$ (otherwise $(F,W)$ is contained in a coherent subsystem $(F,W')$ 
with $\dim W' = 2$ which we have considered already).

Suppose first that $m=0$, so that $G\simeq\cO(a+l)^{n-2}$. Then
$$
e\le(r-1)(a+l)=\frac{r-1}{n-2}d
$$
and the $\alpha$-stability condition for $(F,W)$ is
$$
\frac{e+\alpha}r<\frac{d+2\alpha}n,\ \ \mbox{i.~e.}\ \ (n-2r)\alpha<rd-ne.\eqno(18)
$$
This is true if
$$
(n-2r)(n-2)\alpha<r(n-2)d-n(r-1)d,
$$
i.e. $\alpha<\frac{d}{n-2}$, which holds by (12).

When $m\ge1$, we need a refined version of this argument. We want to prove that $(F,W)$ satisfies the $\alpha$-stability 
condition for
$$
 \alpha < \frac{d}{n-2} - \frac{mn}{2(n-2)}.
$$
By (18) this will hold if
$$
(n-2r)\left(\frac{d}{n-2} - \frac{mn}{2(n-2)}\right)\le rd-ne,
$$
i.e.
$$
2[(n-2)e-(r-1)d]\le m(n-2r).\eqno(19)
$$
We shall prove that (19) always holds for $(E,V)$ general.

Suppose first that $1\le r-1\le m$ and write
$$
e=(a+l+1)(r-1)-f,
$$
so that $f$ is the amount by which $e=\deg K$ falls short of its maximum possible value.
Let $s=\mbox{rk}(K\cap\cO(a+l+1)^m)$. Then it is easy to see that
$$
s\ge r-1-f,\ \ \deg(K\cap\cO(a+l+1)^m)\ge e-(r-1-s)(a+l)=s(a+l)+r-1-f.
$$
This implies that $K$ contains a direct factor $\cO(a+l+1)^{r-1-f}$, which is necessarily a subbundle of $\cO(a+l+1)^m$.
Moreover (19) becomes
$$
2[(n-2)((a+l+1)(r-1)-f)-(r-1)((a+l)(n-2)+m)]\le m(n-2r)
$$
which simplifies to
$$
f\ge r-1-\frac{m}2.
$$
We therefore need to show that $f<r-1-\frac{m}2$ is impossible for $(E,V)$ general.

Now the coherent system $(E,V)$ is given by (1), and therefore corresponds to a $2$-dimensional subspace of
$H^1(G^*)$. If $(F,W)$ exists as above, then there exists a subsheaf $K$ of $G$ as above such that the
$2$-dimensional subspace meets the kernel of the restriction map $H^1(G^*)\to H^1(K^*)$ non-trivially.
Now let $S$ be the closed subset of $H^1(G^*)$ defined by
$$
S = \bigcup \ker(H^1(G^*) \ra H^1(K^*)),
$$
the union being taken over all subsheaves $K$ of $G$ of rank $r-1$ and degree $e$. For general $(E,V)$, 
$(F,W)$ cannot exist if
$$
\mbox{codim}(S,H^1(G^*)) \geq 2.\eqno(20)
$$

One can obtain an estimate of this codimension by computing the dimension of the appropriate Quot-scheme of $G$.
This turns out not to be quite good enough for our purpose. However, from the above discussion, we need only to
prove that (20) holds under the assumption $f<r-1-\frac{m}2$. But now we note that
$$
\ker(H^1(G^*)\ra H^1(K^*))\subset\ker(H^1(G^*)\ra H^1(K'^*))
$$
for some subbundle $K'$ of $\cO(a+l+1)^m$ isomorphic to $\cO(a+l+1)^{r-1-f}$. Since $H^1(G^*)\ra H^1((\cO(a+l+1)^m)^*)$
is surjective, it follows that
$$
\mbox{codim}(S,H^1(G^*)) \geq \mbox{codim}(S',H^1((\cO(a+l+1)^m)^*)),\eqno(21)
$$
where
$$
S' = \bigcup \ker(H^1((\cO(a+l+1)^m)^*) \ra H^1(K'^*)),
$$
the union being taken over all subbundles $K'$ of $\cO(a+l+1)^m$ which are isomorphic to $\cO(a+l+1)^{r-1-f}$.
But this codimension is easy to estimate; in fact, writing $p=r-1-f$,
$$
\mbox{codim}(S',H^1((\cO(a+l+1)^m)^*))\ge\dim H^1(K'^*)-p(m-p)=p(a+l-m+p).\eqno(22)
$$
Now $a\ge l\left(\frac{n-2}2\right)+\frac{m}2$ by (3) and $p=r-1-f>\frac{m}2>0$ by assumption. So
$$
p(a+l-m+p)>p\left(l\left(\frac{n-2}2\right)+\frac{m}2+l-m+\frac{m}2\right)=pl\left(\frac{n}2\right)\ge1.\eqno(23)
$$
The required inequality (20) now follows from (21), (22) and (23). This completes the proof in the case $1\le r-1\le m$.\\

Suppose finally that $m\ge1$ and $r-1>m$. We then write
$$
e=(a+l)(r-1)+m-f
$$
and define $s$ as above. Then
$$
s\ge m-f,\ \ \deg(K\cap\cO(a+l+1)^m)\ge e-(r-1-s)(a+l)=s(a+l)+m-f.
$$
Now (19) becomes
$$
2[(n-2)((a+l)(r-1)+m-f)-(r-1)((a+l)(n-2)+m)]\le m(n-2r)
$$
which simplifies to
$$
f\ge \frac{m}2.
$$
We now argue as before with $p=m-f$ and $f<\frac{m}2$. The proofs of (22) and (23) are unchanged.
This completes the proof of Case 3 and hence of the theorem.  \hfill $\square$\\

{\bf Remark 5.5:} (a) We have already identified one exceptional case $n=4$, $d=6$. In this case (13) holds and $l\ge1$. So the 
proof works at all points other than the one identified above. There do not exist $\alpha$-stable coherent systems of type $(4,6,2)$,
but there do exist $\alpha$-semistable systems provided that the analogue of (12) for semistable coherent systems holds, 
i.~e. $1\le\alpha\le3$.

(b) The condition (13) is used at only two points in the proof. These correspond to two further exceptional cases. The first is
$n=3$, $d=2$ (see the end of Case 1). Here $E\simeq\cO(1)^2\oplus\cO$, which fails several tests for $\alpha$-stability,
in particular Lemma 2.1. In fact, there exists a subsystem $(F,W)$ with $F\simeq\cO(1)^2$, $\dim W=1$, for which the
$\alpha$-semistability condition is $\alpha\ge2$. Since $\alpha\le\frac{d}{n-k}=2$, it follows that $(E,V)$ is
$\alpha$-semistable only for $\alpha=2$.

The second case occurs at the end of Case 2. Here $n=2r$, $E\simeq\cO(r-1)^{2r}$ and $(E,V)$ is $\alpha$-semistable in the range
$0\le\alpha\le r$, but never $\alpha$-stable.\\

To complete this section, we consider the case $n=k=2$.\\

{\bf Proposition 5.6:} {\it $G(\alpha;2,d,2)\ne\emptyset$ if and only if $d>2$ and $\frac{t}2<\alpha$.}\\

{\it Proof:} Suppose first that $G(\alpha;2,d,2)\ne\emptyset$. The Brill-Noether condition (13) gives $d\ge2$. However, if $d=2$,
it is easy to see that $(E,V)$ possesses a subsystem $(F,W)$ with $F$ a line bundle of degree $\ge1$ and $\dim W=1$; this contradicts 
$\alpha$-stability for all $\alpha$. The second condition is just Proposition 4.1.

For the converse, suppose first that $t=0$ and $E\simeq\cO(a)\oplus\cO(a)$ with $a\ge2$. For general $V$ 
and any subsystem $(F,W)$
with $\dim W=1$, (16) and (17) imply that $\deg F\le a-1$. So the $\alpha$-stability condition is $a-1<a$. 

Now suppose that $t=1$ and $E\simeq\cO(a)\oplus\cO(a-1)$ with $a\ge2$. For a general choice of $V$, the map $V\ra H^0(\cO(a-1))$ 
is injective. 
So $F\ne\cO(a)$ and hence $\deg F\le a-1$. Now the $\alpha$-stability condition is $a-1<a-\frac12$. \hfill$\square$\\

\section{Some other cases}

For other values of $k$, we do not have complete results, but there are some cases in which we can obtain partial
information, proving the non-emptiness of some of the $G(\alpha;n,d,k)$. \\

We begin with the case $k=n-1$.\\

{\bf Proposition 6.1:} {\it $G(\alpha;n,d,n-1)$ is non-empty for some $\alpha$ if and only if $d\ge n$. Moreover, in
this case, the upper bound for $\alpha$ is precisely $d$.}\\

{\bf Remark 6.2:} Note that this does not contradict (7) since in this case $m$ is always $0$.\\

{\it Proof}: By \cite[Remark 5.5]{2}, the ``large $\alpha$'' moduli space is always non-empty for $d\ge n$. 
The necessity of the condition $d\ge n$ follows from Corollary 3.3 (indeed $d\ge n$ is precisely the Brill-Noether condition).
\hfill$\square$\\

Next we consider the case $k=n$. We have already solved this case completely when $n=2$ (Proposition 5.6). For $n>2$, 
we have the following partial result.\\

{\bf Proposition 6.3:} {\it  $G(\alpha;n,d,n)$ is non-empty for some $\alpha$ if and only if $d>n$. Moreover, in this
case, there is no upper bound on $\alpha$.}\\

{\it Proof}: Both statements follow from \cite[Theorem 5.6 and Remark 5.7]{2}.\hfill$\square$\\

Finally suppose $k=n+1$.\\

{\bf Proposition 6.4:} {\it $G(\alpha;n,d,n+1)$ is non-empty for some $\alpha$ if and only if $d\ge n$. Moreover, if
$d\ge n$ and we 
write $d=na-t$ with $0\le t<n$, then $G(\alpha;n,d,n+1)$ is always non-empty if $t<\alpha$.}\\

{\bf Remark 6.5:} If $t=0$, the lower bound on $\alpha$ is of course precise. For $t>0$, we know by Propositions 4.1 and 6.4
that the precise lower bound on $\alpha$ lies in the interval $[\frac{t}{n+1},t]$.\\

{\it Proof}: Note that the Brill-Noether condition of Corollary 3.3 is precisely $d\ge n$. The first sentence now follows
directly from \cite[Theorem 5.11]{2}. For the second part, suppose $(F,W)$ is a proper subsystem of $(E,V)$, where $E$ is of
generic splitting type and $V$ globally generates $E$. If $\mbox{rk}F=r$, then $\deg F\le ra$ and $\dim W\le r$; for the latter, 
note that, if $\dim W>r$, then the image of $V$ in $H^0(E/F)$ has dimension $\le n-r$. Since $\deg E/F>0$, this implies that
$E/F$ cannot be globally generated by $V$, which contradicts the choice of $V$. So the stability condition for
$(F,W)$ holds if
$$
a+\alpha<a-\frac{t}{n}+\alpha\frac{n+1}{n},
$$
i.~e. if $\alpha>t$.\hfill$\square$.\\

{\bf Example 6.6:} In the special case $d=n$, we have $\beta(n,n,n+1)=0$. So the moduli space $G(\alpha;n,n,n+1)$ is a single point. 
The corresponding coherent system (which is $\alpha$-stable for all $\alpha>0$) is isomorphic to $(\cO(1)^n,V)$, where $V$ is 
a subspace of $H^0(\cO(1)^n)$  which globally generates $\cO(1)^n$. The subspace $V$ is not unique but all choices of $V$
give isomorphic coherent systems. One can also construct this system using the dual span construction explained in 
\cite[section 5.4]{2}. In fact we have an exact sequence
$$
0\lra F\lra H^0(\cO(n))\otimes\cO\lra\cO(n)\lra0.
$$
Then $(E,V)$ is isomorphic to the coherent system $(F^*,H^0(\cO(n))^*)$ obtained by dualising this sequence.

Note that this example contrasts with the case $g\ge2$, when $G(\alpha;n,n,n+1)=\emptyset$.\\


\begin{thebibliography}{99}
\bibitem[1]{1} M. F. Atiyah: Vector bundles over an elliptic curve. Proc. London Math. Soc. 7 (1957), 414-452.
\bibitem[2]{5} S. B. Bradlow and O. Garc\'{\i}a-Prada: A Hitchin-Kobayashi correspondence for coherent systems on Riemann
surfaces. J. London Math. Soc. 60 (1999), 155-170.
\bibitem[3]{7} S. B. Bradlow and O. Garc\'{\i}a-Prada: An application of coherent systems to a Brill-Noether problem.
J. Reine Angew. Math. 551 (2002), 123-143.
\bibitem[4]{2} S. B. Bradlow, O. Garc\'{\i}a-Prada, V. Mu\~noz and P. E. Newstead: Coherent systems and Brill-Noether theory. 
Internat. J. Math. 14 (2003), 1-50.
\bibitem[5]{6} O. Garc\'{\i}a-Prada: Dimensional reduction of stable bundles, vortices and stable pairs.
Internat. J. Math. 5 (1994), 1-52.
\bibitem[6]{3} S. S. Shatz: On subbundles of vector bundles over $\PP_n$. J. Pure Appl. Algebra 10 (1977), 315-322.
\end{thebibliography}
\end{document}